\documentclass[letter,11pt, reqno]{amsart}
\usepackage[margin=3cm]{geometry}
\usepackage{amsthm}

\usepackage{algorithm}
\usepackage{algpseudocode}
\usepackage{tikz}

\usepackage{epstopdf}% To incorporate .eps illustrations using PDFLaTeX, etc.

\usepackage{xcolor}
\usepackage{url}
\usepackage{hyperref}
\usepackage[caption=false]{subfig}% Support for small, `sub' figures and tables

\usepackage[numbers,sort&compress]{natbib}% Citation support using natbib.sty
\bibpunct[, ]{[}{]}{,}{n}{,}{,}% Citation support using natbib.sty
\theoremstyle{plain}% Theorem-like structures provided by amsthm.sty
\theoremstyle{definition}

\theoremstyle{remark}

\usepackage{adjustbox}

\let\origmaketitle\maketitle
\def\maketitle{
	\begingroup
	\def\uppercasenonmath##1{} % this disables uppercasing title
	\let\MakeUppercase\relax % this disables uppercasing authors
	\origmaketitle
	\endgroup
}

\begin{document}

\title[]{\huge Optimal Embedding of Wiring Diagrams in Constrained Three-Dimensional Spaces}

\author[V. Blanco, G. Gonz\'alez, \MakeLowercase{and} J. Puerto]{
{\large V\'ictor Blanco$^{\dagger}$, Gabriel Gonz\'alez$^{\dagger,\ddagger}$, and Justo Puerto$^{\ddagger}$}\medskip\\
$^\dagger$Institute of Mathematics (IMAG), Universidad de Granada\\
$^\dagger$Institute of Mathematics (IMUS), Universidad de Sevilla\\
\texttt{vblanco@ugr.es}, \texttt{ggdominguez@ugr.es}, \texttt{puerto@us.es}
}

\maketitle

\begin{abstract}
This paper investigates the \emph{Wiring Diagram Problem} (WDP), a three-dimensional layout design problem arising in industrial applications such as cable harness design and pipeline routing in constrained environments. In these settings, hierarchical tree-like systems composed of supply units, intermediate devices (e.g., valves or junctions), and terminal components must be spatially arranged and interconnected while satisfying stringent engineering requirements, including safety separation distances, obstacle avoidance, geometric feasibility, and constructibility constraints. 

We develop an optimization-based framework that formulates the WDP as a mixed-integer linear programming model capturing both topological and spatial design requirements within a unified formulation. To address the combinatorial and geometric complexity of three-dimensional routing, the feasible design space is discretized into structured network graphs that preserve engineering constraints while reducing dimensionality. 

The resulting model minimizes total cable or pipeline length while ensuring compliance with all technical specifications. Computational experiments on representative industrial instances demonstrate the robustness and practical applicability of the proposed approach for automated layout generation.
\end{abstract}


\section{Introduction \label{introduction}}

In many industrial engineering applications, wiring and pipeline layouts constitute a fundamental component of system design. A \emph{wiring diagram} specifies how supply units, intermediate devices (e.g., valves, junctions, or distribution boxes), and terminal elements are interconnected, defining both the topology and the physical routing of cables or pipes within a three-dimensional environment. Beyond representing connectivity, these diagrams must satisfy stringent engineering requirements, including safety separation distances, obstacle avoidance, accessibility for maintenance, and constructibility constraints imposed by industry standards. As industrial systems become increasingly complex, the need for systematic and computationally efficient tools to support layout determination has grown significantly, motivating the development of optimization-based design methodologies.

Designing such layouts is a complex and labor-intensive task. In practice, engineers must route multiple services (e.g., electricity, water, communication lines, or hydraulic systems) through geometrically constrained spaces that may include structural elements, equipment, and restricted zones. The resulting design must simultaneously satisfy functional requirements and minimize installation cost, typically measured by total cable or pipeline length, while preserving feasibility and maintainability. These challenges are particularly relevant in domains such as naval architecture, wind farm infrastructure, oil and gas transportation, and large-scale industrial facilities, where spatial congestion and strict safety regulations significantly increase design complexity~\cite{blanco2021network,blanco2023pipelines,fischetti2018optimizing,liu2024improved,soliman1986gas,wang2022pipeline}. Similarly, in chemical and process engineering, the joint optimization of equipment placement and pipe routing under safety separation constraints has been formulated as mixed-integer programming problems~\citep{guirardello2005,georgiadis1999}, further illustrating the broad industrial relevance of constrained layout and routing optimization. In related contexts, cost-optimal structural and modular system design problems have also been addressed through mathematical optimization techniques, such as geometric programming approaches for engineering layout optimization~\cite{paul1972geometric}.

From a methodological perspective, routing problems of this type are often related to shortest-path or Steiner tree formulations on graphs, extended to incorporate engineering constraints such as bend limitations, collision avoidance, or service separation requirements.  Tree-structured network layout problems have been investigated using both evolutionary and exact optimization methods. For instance, \cite{walters1987evolutionary} proposed evolutionary strategies for optimal tree network layout, while large-scale gas pipeline design problems have been addressed via mathematical optimization  formulations~\cite{soliman1986gas}. More recently, spatially flexible tree embedding models with placement regions have been studied under neighborhood-based Steiner formulations~\cite{blanco2024ftstn}. 
A classical tool for discretizing rectilinear Steiner problems is the Hanan grid~\citep{hanan1966steiner}, whose optimality properties have been extended to broader orientation metrics~\citep{mullerhannemann2015}. A comprehensive survey of automatic pipe routing methods, covering discretization strategies, Steiner tree formulations, obstacle avoidance, and safety constraints, is provided in~\citep{APRsurvey2023}.

Recent contributions have addressed related design problems across several engineering domains using mathematical optimization and metaheuristic approaches. The rapid progress in optimization algorithms, together with the increasing computational power of modern solvers and hardware, has considerably expanded the range and scale of engineering layout problems that can now be treated using optimization-based decision-support tools. In naval engineering, pipe routing optimization has been tackled using pathfinding algorithms such as Dijkstra's method~\citep{parkstorch2002,gunawan2022}, expert systems combined with multi-objective evolutionary search~\citep{lee2023shippiping}, and ant colony optimization for branch pipe routing~\citep{jiang2015}. Network-flow-based MILP models for multicommodity pipeline routing in ship compartments have also been recently developed~\citep{blanco2021network,blanco2023pipelines}. In the cable harness design domain, automated cable harness routing in three-dimensional environments has been addressed through graph-based multi-objective formulations combining shortest-path and Steiner tree models~\citep{karlsson2024}, and collision-aware routing strategies using improved pathfinding algorithms with B-spline optimization~\citep{kim2024}, as well as multi-objective particle swarm optimization for multi-branch cable harness layout in aircraft~\citep{zhang2021cablebranch}. Automatic planning for robotic wiring execution based on digital design documents has also been investigated~\cite{shneor2023roboticwiring}. In the wind energy sector, mixed-integer linear programming (MILP) formulations have been extensively applied to wind farm cable routing, including models with obstacle avoidance via Steiner points and power loss optimization~\citep{wedzik2016,fischetti2018optimizing}, as well as integrated frameworks that simultaneously determine turbine placement and cable routing decisions~\citep{fischetti2022}. Notably, MILP-based cable layout models with engineering constraints have been shown to be effective in this journal~\citep{klein01022021}. In chemical and process plant design, three-dimensional pipe routing with safety and flexibility constraints has been studied using MILP~\citep{guirardello2005}, constraint programming~\citep{belov2017piprouting}, and hybrid LP and MILP formulations~\citep{stanczak2020piprouting}. While these works demonstrate the growing potential of optimization-driven automation in layout design, they typically focus on specific routing tasks, individual pipe or cable connections, robotic execution planning, or simplified geometric settings. The systematic embedding of hierarchical wiring diagrams within fully constrained three-dimensional industrial spaces, explicitly enforcing safety separation and constructibility requirements within a unified optimization framework, remains largely unexplored.

This paper addresses the \emph{Wiring Diagram Problem} (WDP), which consists of embedding a hierarchical, tree-structured interconnection scheme into a constrained three-dimensional design domain. The problem integrates topological requirements, imposed by the predefined wiring structure, with spatial and engineering constraints that govern feasible placement and routing. In particular, intermediate components must be located within admissible regions, connections must avoid obstacles and respect separation distances, and the resulting layout must remain constructible according to industrial practice.

To tackle this problem, we develop an optimization-based framework that combines structured spatial discretizationn with mixed-integer linear optimization. The three-dimensional design domain is discretized into tailored network graphs that encode feasible routing corridors and admissible locations for intermediate components. This discretization significantly reduces geometric complexity while preserving engineering feasibility and enabling tractable optimization models. In this network representation, we formulate a MILP model that enforces hierarchical connectivity, spatial compliance, routing feasibility, and safety-distance constraints within a unified mathematical optimization framework.

The proposed model minimizes the total length of the cable or pipeline subject to all technical constraints, thereby supporting the generation of feasible and cost-efficient layout. The approach leverages modern mathematical optimization tools to produce layouts that satisfy engineering requirements while improving design efficiency compared with manual or heuristic procedures. In this way, the methodology provides a systematic and reproducible decision-support tool for practical engineering layout design.

The remainder of the paper is organized as follows. Section~\ref{sec:spp} formally defines the Wiring Diagram Problem. Section \ref{sec:discretization} presents the discretization based on graphs of the design space, and Section~\ref{modelo} introduces the proposed optimization model. Sections~\ref{sec:computational} and~\ref{sec:case} report computational experiments and an industrial case study, illustrating the capability of the approach to generate feasible and efficient layouts in representative industrial scenarios.

\section{Wiring Diagram Problem}\label{sec:spp}

In this section, we introduce the problem under analysis and establish the notation used throughout the paper.

\begin{figure}[h]
\begin{center}
    \includegraphics[width=\textwidth]{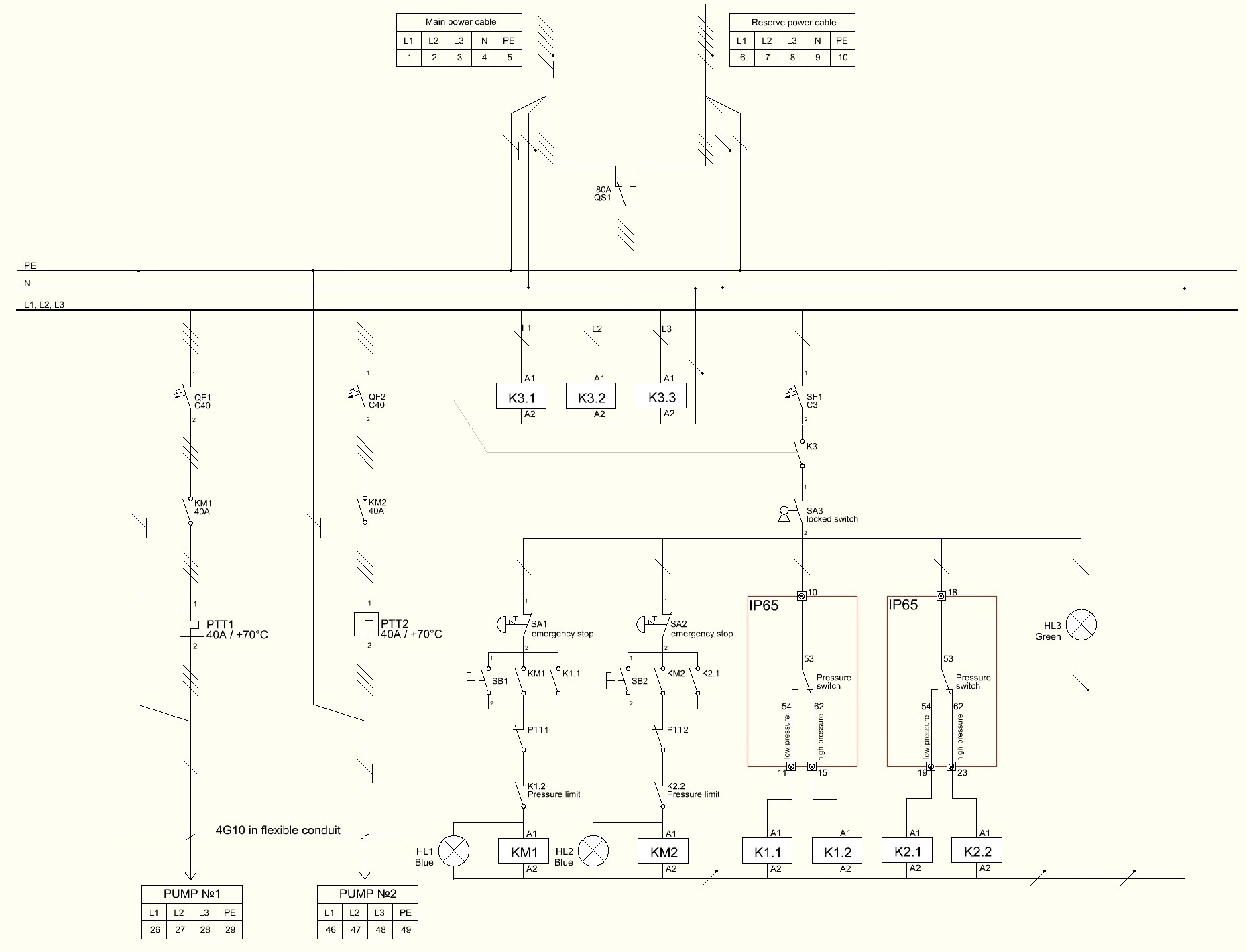}
\end{center}
\caption{Example of wiring diagram (source: \url{https://commons.wikimedia.org/wiki/Motor\_control\_centre\#Wiring\_diagrams}).\label{fig:wd0}}
\end{figure}

The Wiring Diagram Problem (WDP) considered in this work consists of embedding a predefined hierarchical interconnection scheme into a constrained three-dimensional engineering domain. 
A wiring diagram is a hierarchical reference structure that specifies the connectivity between supply sources, intermediate devices (e.g., valves or junction systems), and terminal components in an industrial installation. An example of a wiring diagram is shown in Figure \ref{fig:wd0}. This figure depicts a motor control center in a pump station, where a main power cable distributes electricity to several demand points through intermediate devices arranged according to a prescribed hierarchy.

Beyond connectivity, a practical implementation must satisfy stringent engineering requirements, including spatial feasibility, safety separation distances, obstacle avoidance, and constructibility constraints. Therefore, the problem is not purely topological, but inherently geometric and constrained.

To describe the problem in a unified mathematical framework, we represent the wiring diagram as a finite collection of rooted directed trees, that is, a \emph{forest}. Each tree encodes a hierarchical branch of the installation. The root of each tree corresponds to a supply source connected to an existing main pipeline, while the leaves represent terminal devices with fixed spatial coordinates. Intermediate nodes correspond to valves or junction components that must be located within admissible spatial regions.

More precisely:
\begin{itemize}
\item {\bf Root node:} The supply source must be located along a pre-existing main pipeline. These pipelines are assumed to be already designed (e.g., via established methodologies such as \citep{blanco2021network}) and serve as anchoring structures for the branches.
\item {\bf Terminal nodes:} These correspond to demand points whose spatial coordinates are fixed. The design must ensure proper service delivery to each of these locations.
\item {\bf Intermediate nodes:} These represent components such as valves or distribution elements. Each intermediate node is associated with a predefined admissible three-dimensional region within which its exact position must be determined.
\end{itemize}

% Gabriel cambios:
% \L, \C, \I -> \mathcal{}
% \N, \E -> N, E
We denote by $\mathcal{C}$ the set of main pipelines and by $\mathcal{F}$ the forest of rooted trees. Each tree $T \in \mathcal{F}$ is associated with a pipeline $\phi(T) \in \mathcal{C}$. The node set of $T$ is $N(T)$ and its arc set is $E(T)$. The root is denoted by $r(T)$, the set of leaves by $\mathcal{L}(T)$, and the set of intermediate nodes by $\mathcal{I}(T)$. For each non-leaf node $s \in N(T)$, we denote by $\varrho^{-1}(s)$ the set of children of $s$ in $T$.

The WDP therefore integrates a predefined hierarchical topology with spatial placement decisions and routing decisions in three dimensions.

\begin{figure}[h]
\begin{center}
\includegraphics[trim={11cm 1cm 9cm 1cm},clip,width=0.36\textwidth]{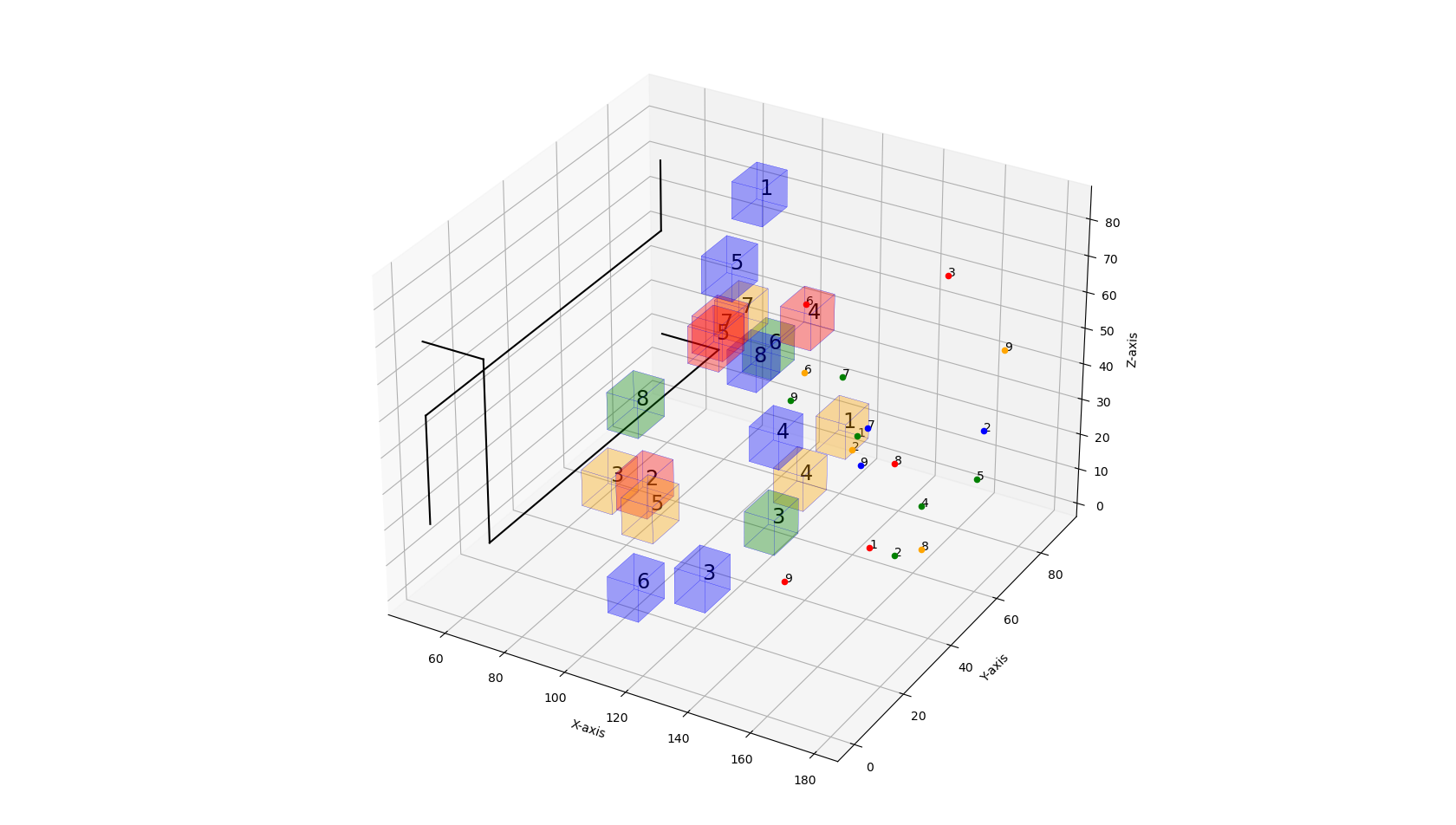}
\includegraphics[trim={2.1cm 0.8cm 2.1cm 0.2cm},clip,width=0.63\textwidth]{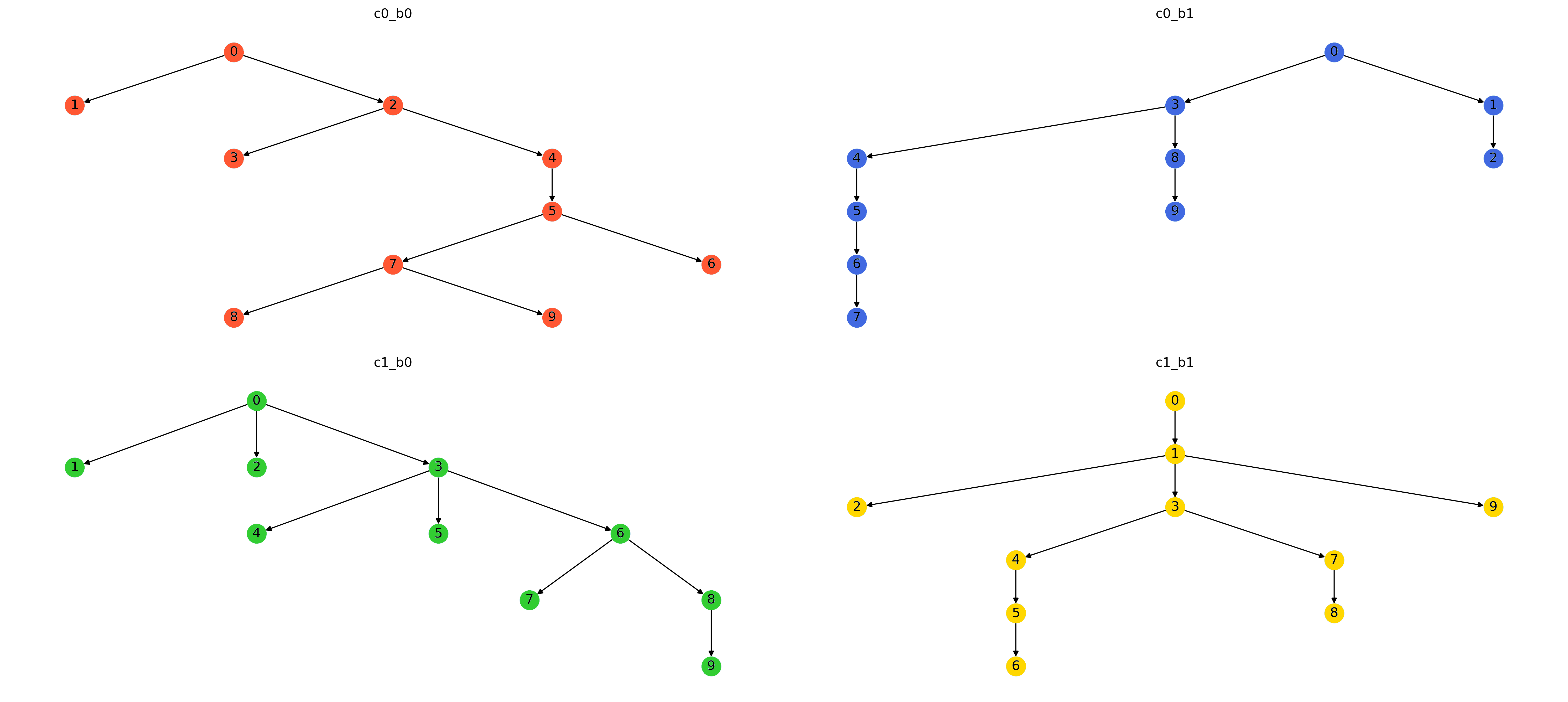}
\caption{\label{fig:input_example} On the left, admissible regions for node placement; on the right, the associated forest structure.}
\end{center}
\end{figure}

Figure \ref{fig:input_example} illustrates an instance with a tree ensemble consisting of four trees associated with two pipelines.

In addition to topological consistency, the following engineering conditions must be satisfied:

\begin{itemize}
\item[-] All nodes must be located within their admissible regions: root nodes along their corresponding pipelines, intermediate nodes within predefined spatial regions, and leaves at fixed coordinates.
\item[-] A minimum safety distance must be maintained between branches and between branches and pipelines, ensuring accessibility and preventing interference.
\item[-] The three-dimensional region may contain obstacles (e.g., structural elements or reserved areas) that must be avoided by all routed connections.
\item[-] When bends are introduced, minimum separation distances between consecutive elbows must be enforced to ensure structural integrity. In the proposed discretization (Section \ref{sec:discretization}), this requirement is implicitly satisfied by the grid resolution, which guarantees a minimum spacing between consecutive direction changes.
\end{itemize}

The objective of the WDP is to determine the spatial placement of all intermediate nodes and the routing of all connections so that the complete forest is embedded in the three-dimensional region while satisfying all engineering constraints and minimizing total cable or pipeline length.

\begin{figure}
\begin{center}
\includegraphics[trim={11cm 1cm 9cm 1cm},clip,width=0.49\textwidth]{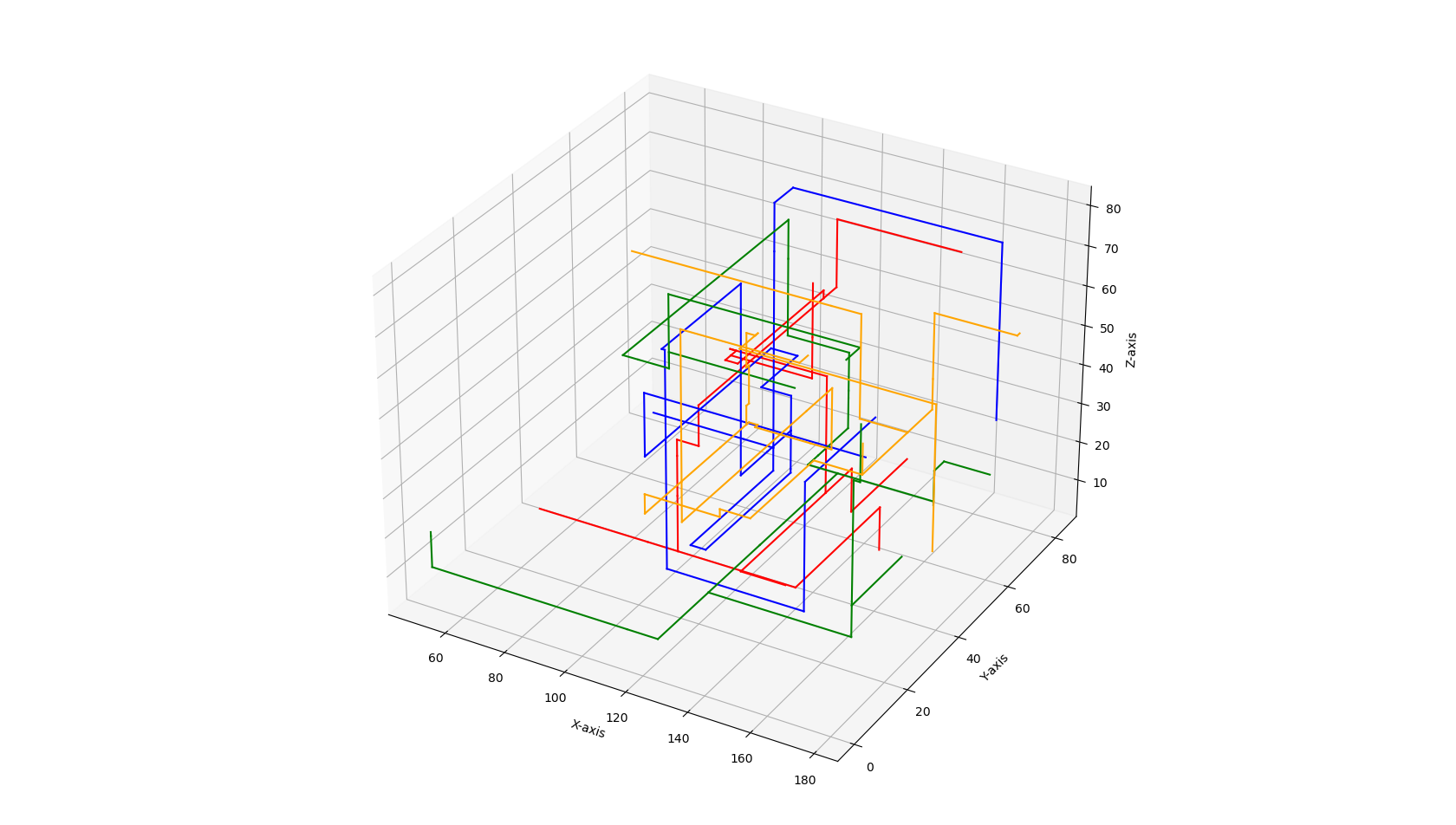}
\includegraphics[trim={11cm 1cm 9cm 1cm},clip,width=0.49\textwidth]{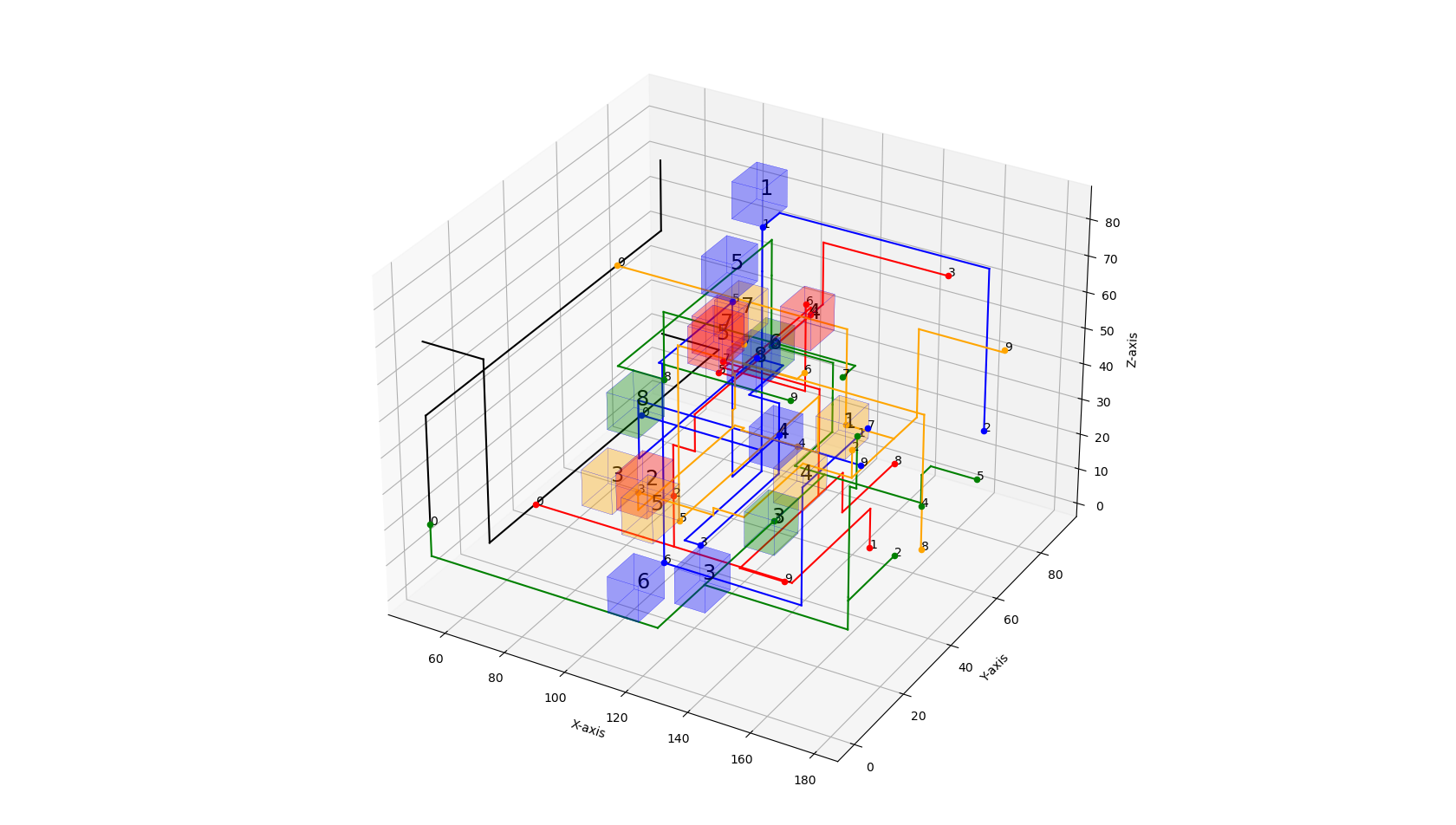}
\caption{\label{solucioncita} Illustration of a feasible solution.}
\end{center}
\end{figure}

Figure \ref{solucioncita} illustrates a feasible layout for the instance in Figure \ref{fig:input_example}.

\section{Graph-based discretization of the design space}\label{sec:discretization}

To transform the continuous three-dimensional layout problem into a tractable optimization model, we discretize the design domain and construct a structured network representation that encodes feasible routing corridors.

Let $\mathcal{R} \subset \mathbb{R}^3$ denote the design region. For each tree $T \in \mathcal{F}$ associated with a pipeline $c=\phi(T)$, we construct an orthogonal grid over $\mathcal{R}$ that includes:
\begin{itemize}
\item extremities and direction-change points of the main pipeline,
\item vertices of admissible regions for intermediate nodes,
\item fixed positions of terminal nodes.
\end{itemize}

This grid induces an undirected graph $\mathcal{G}(T)=(\tilde{N}(T),\tilde{E}(T))$. Nodes and edges intersecting obstacles are removed to preserve feasibility. The global discretized structure is $\mathcal{G}=\bigcup_{T\in\mathcal{F}}\mathcal{G}(T)$. Each edge $e\in\tilde{E}(T)$ is associated with a length $d_e$. For each intermediate node $s \in \mathcal{I}(T)$, we denote by $\mathcal{R}_s \subset \tilde{N}(T)$ the set of grid vertices lying within the admissible region of $s$.

A naive discretization would generate extremely large graphs, rendering direct optimization impractical. To mitigate this difficulty, we exploit the hierarchical structure of the wiring diagram to construct reduced grids tailored to each parent-children configuration.

For a node $s$ and its children $\varrho^{-1}(s)$ in tree $T$, we build a three-dimensional Hanan grid~\citep{hanan1966steiner} generated by the admissible regions of $s$ and its children. Under $\ell_1$-based distances, optimal rectilinear connections lie within this grid. 

Unlike classical independent Steiner tree constructions, however, branches in the WDP compete for shared space and must jointly satisfy safety constraints. Therefore, the discretization serves as a structured restriction of the continuous domain that reduces dimensionality while preserving practical feasibility. Although global optimality in the continuous domain cannot be guaranteed, the construction provides a tractable and engineering-consistent search space.

\begin{figure}
\begin{center}
\includegraphics[trim={4.82cm 0.1cm 6.53cm 0.15cm},clip,width=0.28\textwidth]{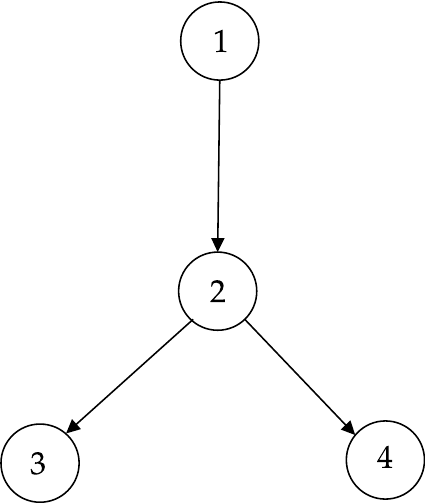}
\includegraphics[trim={2.21cm 0.1cm 4.19cm 0.15cm},clip,width=0.3\textwidth]{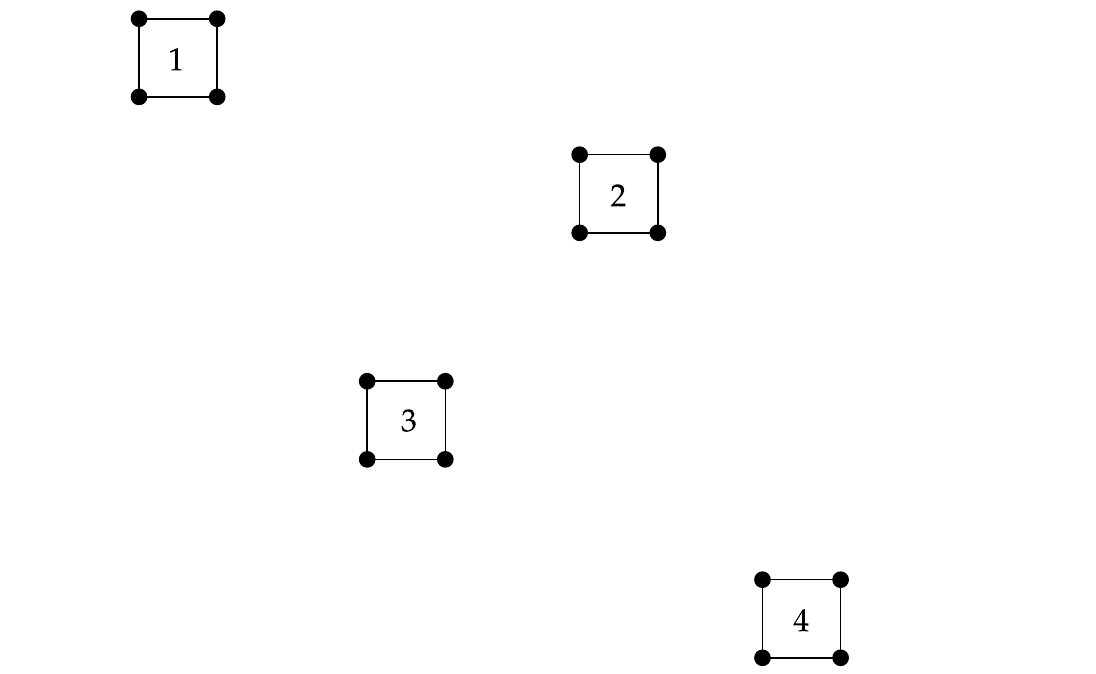}
\includegraphics[trim={2.74cm 0.1cm 3.65cm 0.05cm},clip,width=0.3\textwidth]{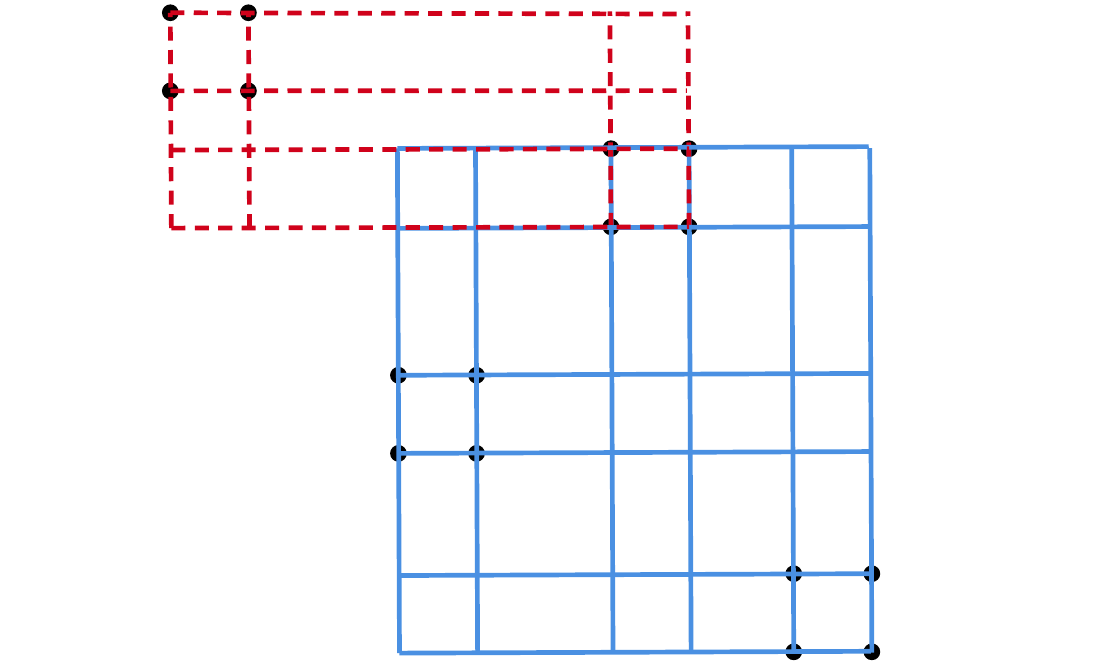}
\caption{Example of solution space discretization.}
\label{crear_mallado}
\end{center}
\end{figure}

Figure~\ref{crear_mallado} illustrates the discretization process for a subtree.

\begin{figure}
\begin{center}
\includegraphics[trim={0cm 0.1cm 2.4cm 0.1cm},clip,width=0.7\textwidth]{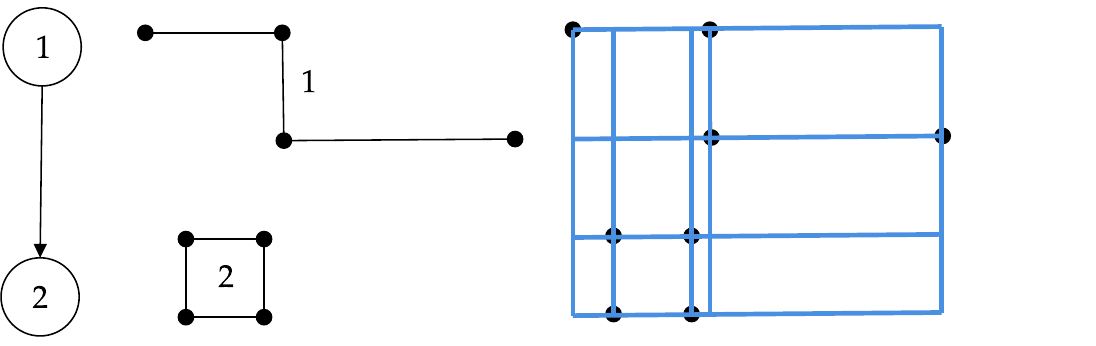}
\caption{Discretization at a root node.}
\label{colector_mallado}
\end{center}
\end{figure}

\begin{figure}
\begin{center}
\includegraphics[trim={0cm 0.1cm 4.05cm 0cm},clip,width=0.7\textwidth]{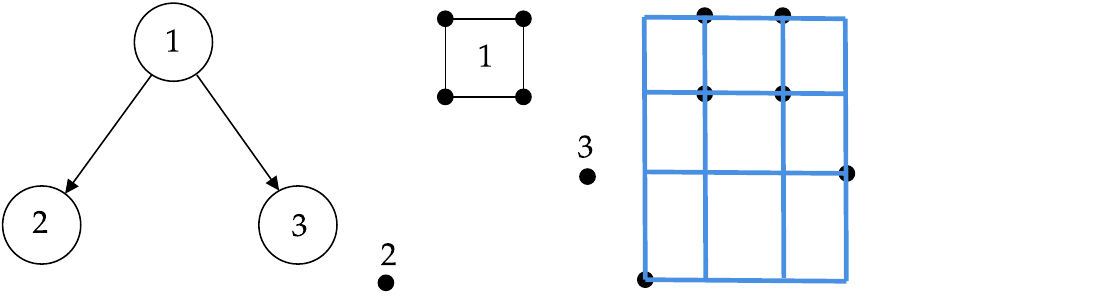}
\caption{Discretization at a leaf node.}
\label{hojas_mallado}
\end{center}
\end{figure}

Special cases arise when the parent is a root node (pipeline extremities are used instead of valve regions) or when a child is a leaf (a fixed point replaces an admissible region), as illustrated in Figures~\ref{colector_mallado} and~\ref{hojas_mallado}.

The resulting collection of structured grids defines a finite network on which the mixed-integer optimization model is formulated. This discretization enables the integration of routing, placement, and safety constraints within a unified and computationally tractable framework.

\section{Mathematical Optimization Model\label{modelo}}

The Wiring Diagram Problem is formulated as a mixed-integer linear programming (MILP) model defined over the reduced discretization graph introduced in the previous section. The model integrates three tightly coupled design decisions: (i) spatial placement of intermediate components, (ii) routing of connections between hierarchical nodes, and (iii) enforcement of safety and constructibility constraints.

The routing metric is based on the $\ell_{1}$-norm, which is consistent with rectilinear industrial layouts and compatible with the orthogonal grid structure used in the discretization.

While the original discretization of the three-dimensional region may generate very large graphs, the structured reduction strategy described in the previous section enables the construction of a tractable network on which an exact optimization model can be solved efficiently.

To describe safety interactions between routed segments, we define:

\begin{itemize}
\item[] $\delta_{aa'}$: the $\ell_1$-distance between arcs $a$ and $a'$ in $\mathbb{R}^3$.
\item[] $\Delta$: the prescribed minimum safety distance between any two routed segments.
\end{itemize}
We also use $\delta^-(v)$ and $\delta^+(v)$ to denote, respectively, the sets of incoming and outgoing arcs at vertex $v$ in the discretized graph.

\subsection*{Decision variables}

The model uses the following binary decision variables:

\begin{align*} 
f_{a}^{st} &= 
\begin{cases}
1 & \text{if the path from $s$ to $t$ uses arc $a$},\\
0 & \text{otherwise}
\end{cases}
&\forall a \in \tilde{E}(T),\ s,t \in \tilde{N}(T),\ T\in\mathcal{F},\\[6pt]
y_{v}^{s} &= 
\begin{cases}
1 & \text{if intermediate node $s$ is placed at vertex $v$},\\
0 & \text{otherwise}
\end{cases}
&\forall v \in \mathcal{R}_s,\ s \in \mathcal{I}(T),\\[6pt]
x_{a} &= 
\begin{cases}
1 & \text{if arc $a$ is installed in the final layout},\\
0 & \text{otherwise}
\end{cases}
&\forall a \in \tilde{E}(T),\ T\in\mathcal{F}.
\end{align*}

Variable $y$ determines the physical placement of valves or intermediate components. Variable $x$ identifies which routing segments are constructed. Variable $f$ ensures connectivity between hierarchical elements by enforcing feasible paths within the selected network. The flow variables $f_a^{st}$ are defined for each parent--child pair $(s,t) \in E(T)$ in the tree, representing the unit flow that must be routed from the selected position of $s$ to the selected position of $t$.

\subsection*{Constraints}
\subsubsection*{Component placement constraints}
Each intermediate component must be assigned to exactly one admissible spatial location:
\begin{align}
\sum_{v \in \mathcal{R}_{s}} y_{v}^{s} = 1,
\qquad \forall s \in \mathcal{I}(T),\ T \in \mathcal{F}.
\end{align}
This guarantees a unique physical position for every valve or intermediate device within its designated region.

\subsubsection*{Routing and connectivity constraints}

Connectivity between parent and child nodes is enforced using a single-commodity flow representation on the discretized graph.

Flow conservation is imposed as:
\begin{align}
\sum_{a \in \delta^-(v)} f^{st}_{a} -
\sum_{a \in \delta^+(v)} f^{st}_{a}
&= y_{v}^{s},
&&\forall s,t,\ v \in \mathcal{R}_{s}, \\
\sum_{a \in \delta^-(v)} f^{st}_{a} -
\sum_{a \in \delta^+(v)} f^{st}_{a}
&= -y_{v}^{t},
&&\forall s,t,\ v \in \mathcal{R}_{t}, \\
\sum_{a \in \delta^-(v)} f^{st}_{a} -
\sum_{a \in \delta^+(v)} f^{st}_{a}
&= 0,
&&\text{otherwise}.
\end{align}
These constraints ensure that exactly one unit of flow is sent from the selected position of node $s$ to the selected position of node $t$, thereby enforcing hierarchical connectivity.

Routing can only occur on installed arcs:
\begin{align}
f^{st}_{a} \le x_a,
\qquad \forall a,s,t.
\end{align}
Thus, only physically constructed segments can carry connections. For leaf nodes:
\begin{align}
\sum_{a \in \delta^-(\ell)} x_a = 1,
\qquad
\sum_{a \in \delta^+(\ell)} x_a = 0,
\qquad \forall \ell \in \mathcal{L}(T).
\end{align}
These constraints enforce that each terminal device is connected by exactly one incoming arc and does not generate outgoing branches.
\subsubsection*{Unidirectional flow constraints}
To prevent flow from traversing both orientations of the same edge simultaneously, we impose:
\begin{align}
f^{st}_{a} + f^{st}_{a'} \le 1,
\qquad \forall a,s,t.
\end{align}
\subsubsection*{Safety constraints}

Minimum separation between routed segments is enforced through pairwise exclusion constraints:
\begin{align}
\sum_{a' : \delta_{aa'} < \Delta} x_{a'} \le M(1-x_a),
\qquad \forall a.
\end{align}
where $M$ is a sufficiently large constant (e.g., the total number of arcs in $\tilde{E}$). If arc $a$ is selected, no other arc within the prescribed safety distance can be activated simultaneously. This guarantees compliance with maintenance and interference regulations.

Similar constraints prevent routed branches from violating safety distances with respect to main pipelines.
\begin{align}
\sum_{c' \in \mathcal{C}:\, c' \neq c}\;\sum_{T' \in \mathcal{F}:\, \phi(T')=c'}\;\sum_{\substack{a \in \tilde{E}(T'):\\ \Delta^{c}_{a} < \Delta}} x_{a} = 0, \qquad \forall c \in \mathcal{C},
\end{align}
where $\Delta^{c}_{a}$ denotes the $\ell_1$-distance from arc $a$ to pipeline $c$.

In practice, the pairwise safety constraints grow quadratically with the number of arcs. To avoid enumerating all arc conflicts a priori, constraints (8) and (9) are incorporated dynamically via lazy constraint callbacks during the branch-and-bound process: at each integer-feasible node, safety violations are checked and the corresponding constraints are added only as needed.

\subsection*{Objective function}

The objective is to minimize total installation cost, represented by the total routed length, i.e.:
$$
\sum_{a \in \tilde{E}} d_a x_a.
$$
This objective reflects the primary engineering cost driver: total cable or pipeline length, while all safety and constructibility requirements are strictly enforced through constraints.

The resulting MILP simultaneously determines optimal component placement and routing decisions within a unified framework, providing an automated and engineering-consistent solution to the Wiring Diagram Problem.

\section{Computational Experiments}\label{sec:computational}

In this section, we report computational results on synthetic instances to validate the proposed approach and assess its scalability. In the next section, we apply it to a real industrial case study. The experiments are designed to (i) validate the proposed discretization--MILP framework, (ii) assess scalability with respect to problem size and safety requirements, and (iii) identify practical limitations for highly congested layouts.

All experiments were conducted on a laptop equipped with an 11th Gen Intel Core i7-1165G7 processor (2.80 GHz, 4 cores) and 16 GB of RAM, under Microsoft Windows 11. The MILP models were solved using Gurobi Optimizer version 9.5.0, with a time limit of 3600 seconds and the optimality gap tolerance set to 0\%.

Each instance is defined on a three-dimensional cubic region $\mathcal{R}\subset\mathbb{R}^3$ and contains one or two main pipelines (collectors) located on one side of the cube. From these pipelines, multiple wiring branches must be routed to terminal devices located on the opposite side of the domain. Along each branch, a sequence of intermediate components (e.g., valves/junctions) must be placed within prescribed admissible regions (rectangular boxes) and connected according to a given rooted tree structure. Examples of generated instances are illustrated in Figure~\ref{fig:input_example}.

Main pipelines are generated on an orthogonal grid whose resolution is selected to enforce a minimum bend-to-bend spacing of 10 units. Specifically, each pipeline starts at $(50,0,z)$ and ends at $(50,90,z)$, with $z\in\{0,10,20,30,40,50,60,70,80,90\}$. The pipeline routing is obtained by applying the network-flow methodology proposed in \citep{blanco2021network}, which determines minimum-cost rectilinear routes for each pipeline on a weighted grid while enforcing a minimum separation distance of 6 units between pipelines.

The remaining instance parameters control: (i) the cube dimensions, (ii) the number of branches per pipeline, (iii) the number of nodes per branch (tree size), (iv) the size of admissible regions for intermediate nodes, and (v) the approximate locations of intermediate regions and terminal devices.

For each configuration, we generated 10 independent instances (using different random seeds). For each instance, we solved the resulting MILP on the reduced discretization graph described in Section~\ref{sec:spp}, with a time limit of 3600 seconds. We report: (i) the average solution time (seconds) over the 10 instances, and (ii) the number of instances for which a feasible solution was found within the time limit (denoted \#Feas). 
All computational times correspond to wall-clock seconds.

Table~\ref{tab:results} summarizes the impact of the number of pipelines ($\#c$), the number of branches per pipeline ($\#b$), the number of nodes per branch ($\#n$), and the required safety distance ($\Delta$) on computational performance.

\begin{table}[t]
\centering
\label{tab:results}
\adjustbox{width=0.8\textwidth}{
\begin{tabular}{l l l | r r | r r | r r}
\hline
&&& \multicolumn{6}{c}{Safety distance ($\Delta$)}\\
\textbf{\#c} & \textbf{\#b} & \textbf{\#n} 
& \multicolumn{2}{c|}{\textbf{1}} 
& \multicolumn{2}{c|}{\textbf{3}} 
& \multicolumn{2}{c}{\textbf{5}} \\
& & & Time & \#Feas & Time & \#Feas & Time & \#Feas \\
\hline
\textbf{1} & \textbf{1} & 3  & 0.01 & 10 & 0.01 & 10 & 0.01 & 10 \\
           &            & 5  & 0.05 & 10 & 0.06 & 10 & 0.05 & 10 \\
           &            & 10 & 0.28 & 10 & 0.28 & 10 & 0.28 & 10 \\
           &            & 15 & 0.53 & 10 & 0.53 & 10 & 0.60 & 10 \\
           & \textbf{3} & 3  & 0.05 & 10 & 0.05 & 10 & 0.06 & 10 \\
           &            & 5  & 0.15 & 10 & 0.22 & 10 & 0.21 & 10 \\
           &            & 10 & 7.77 & 10 & 3.58 & 10 & 27.22& 10 \\
           &            & 15 & 3.32 & 10 & 3.80 & 10 & 34.56& 10 \\
           & \textbf{5} & 3  & 0.13 & 10 & 0.13 & 10 & 0.19 & 10 \\
           &            & 5  & 0.36 & 10 & 0.50 & 10 & 0.71 & 10 \\
           &            & 10 & 2.56 & 10 & 5.71 & 10 & 34.49& 10 \\
           &            & 15 & 14.50& 10 & 177.44& 10 & 2711.66& 10 \\\hline
\textbf{2} & \textbf{1} & 3  & 0.03 & 10 & 0.03 & 10 & 0.28 & 10 \\
           &            & 5  & 0.11 & 10 & 0.12 & 10 & 0.13 & 10 \\
           &            & 10 & 0.53 & 10 & 0.50 & 10 & 0.65 & 10 \\
           &            & 15 & 1.96 & 10 & 2.19 & 10 & 2.95 & 10 \\
           & \textbf{3} & 3  & 0.10 & 10 & 0.10 & 10 & 0.14 & 10 \\
           &            & 5  & 0.46 & 10 & 0.60 & 10 & 1.05 & 10 \\
           &            & 10 & 4.23 & 10 & 14.57& 10 & 518.37& 10 \\
           &            & 15 & 10.30& 10 & 228.43& 10 & 2885.18& 10 \\
           & \textbf{5} & 3  & 0.25 & 10 & 0.32 & 10 & 0.44 & 10 \\
           &            & 5  & 1.67 & 10 & 2.98 & 10 & 9.13 & 10 \\
           &            & 10 & 11.91& 10 & 257.72& 10 & 3600.27& 10 \\
           &            & 15 & 206.04&10 & 3600.00& 0 & 3600.00& 0 \\
\hline
\end{tabular}}
\caption{Computational performance as a function of the number of pipelines ($\#c$), branches per pipeline ($\#b$), nodes per branch ($\#n$), and safety distance ($\Delta$). ``Time'' is the average solution time (seconds) over 10 instances. ``\#Feas'' is the number of instances for which a feasible solution was found within 3600 seconds.}
\end{table}

Several conclusions can be drawn from the computational results:
\begin{enumerate}
\item {\it Impact of safety distance:} 
As expected, larger safety distances increase the combinatorial difficulty of the problem because fewer routing segments can be simultaneously activated without violating separation requirements. This effect becomes pronounced when the layout is dense (many branches and nodes), as the feasible routing corridors narrow substantially.
\item {\it Impact of problem size:} 
For small and medium configurations, the model finds feasible solutions quickly across all tested safety distances. As the number of branches ($\#b$) and nodes per branch ($\#n$) increases, the MILP becomes more challenging due to the growth in both routing variables and safety interaction constraints. The most difficult configurations correspond to two pipelines with many branches and large safety distances, where feasible designs may not exist within the discretized space or may require substantial computational effort to identify. All unsolved instances share the same structural profile: maximum complexity ($\#c = 2$, $\#b = 5$, $\#n = 15$) combined with moderate to strict safety distances ($\Delta \geq 3$), confirming that the interaction between spatial conflicts and topological complexity is the primary driver of computational difficulty.
\end{enumerate}

\begin{figure}[t]
    \centering
    \includegraphics[width=\linewidth]{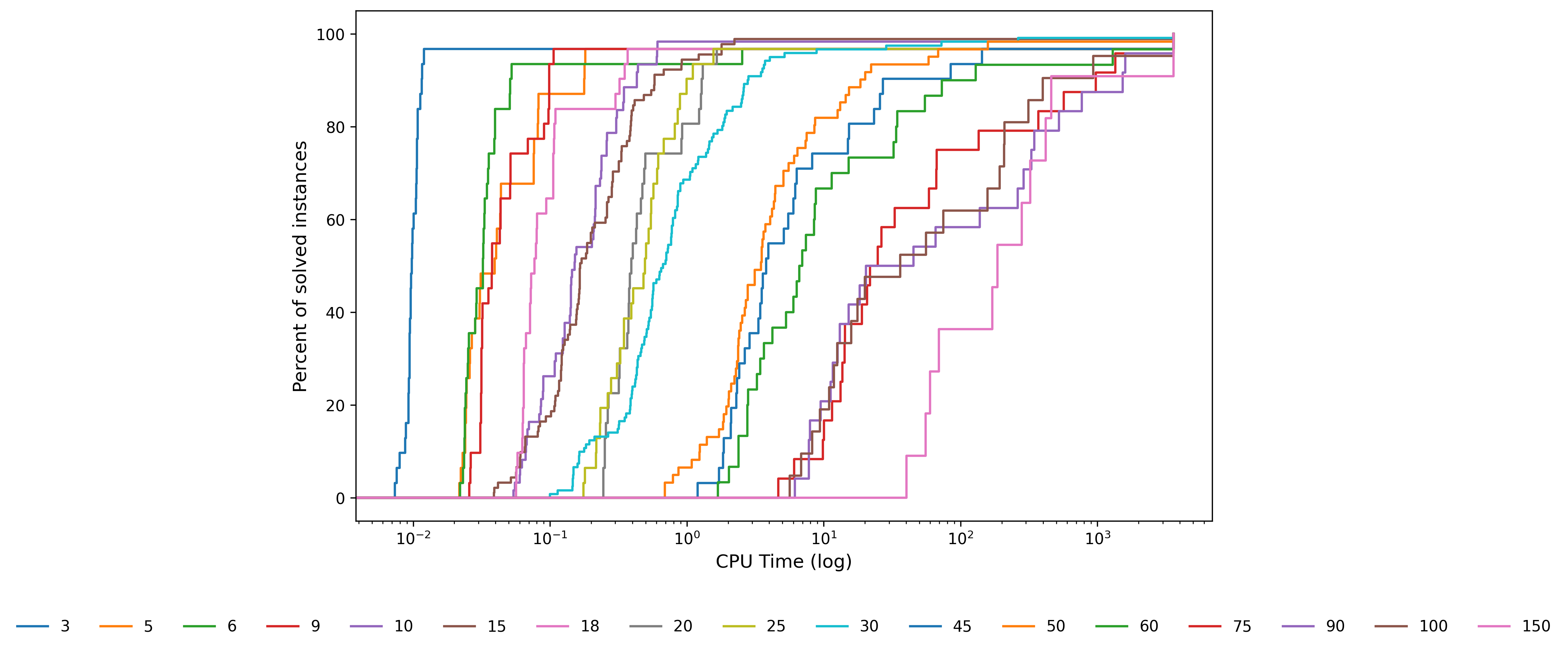}
    \caption{Solution time as a function of the number of pipelines, branches per pipeline, and nodes per branch.
    \label{fig:pp_colectoresxramalxnodo}}
\end{figure}

\begin{figure}[t]
    \centering
    \includegraphics[width=\linewidth]{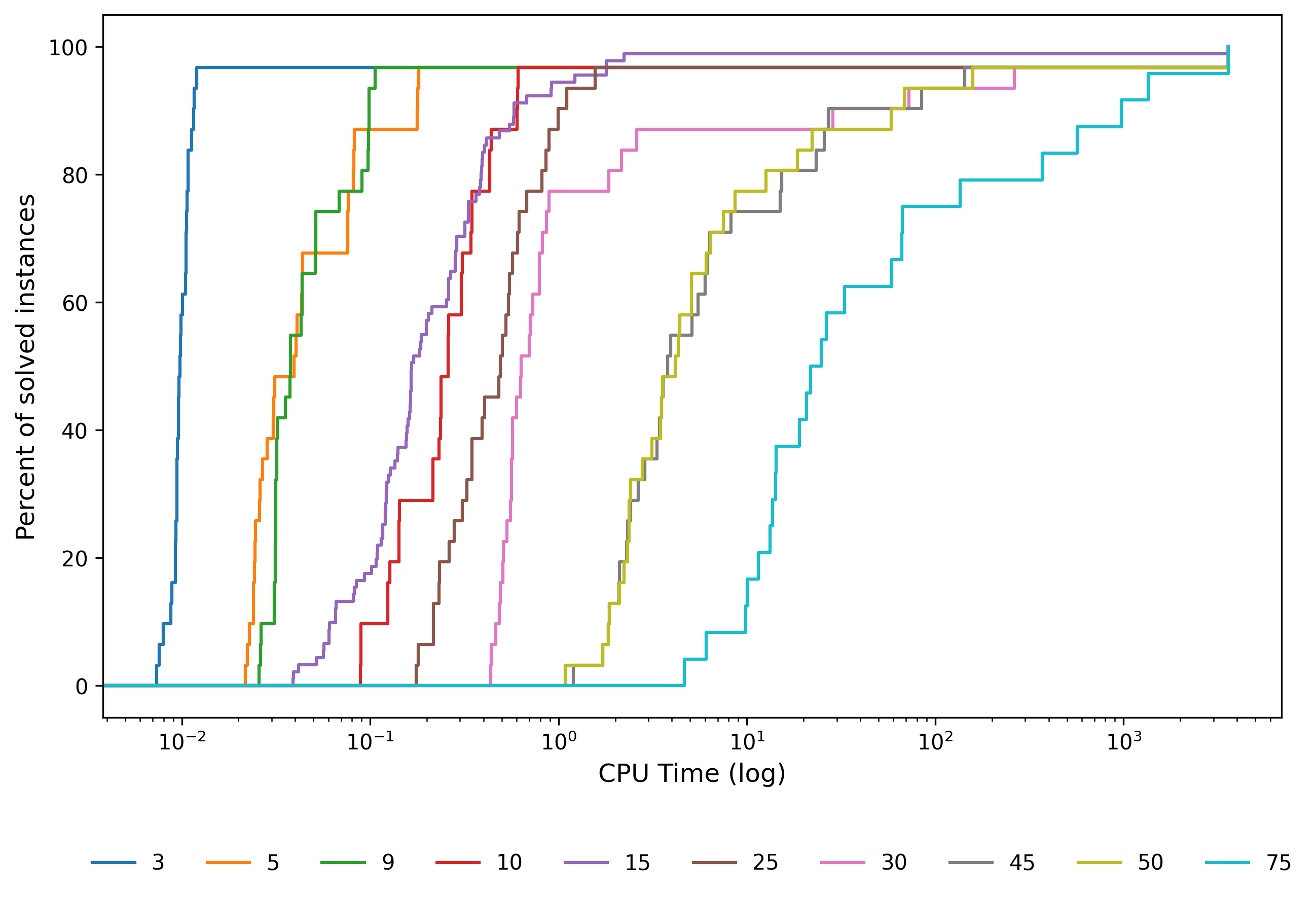}
    \caption{Solution time as a function of branches per pipeline and nodes per branch (one pipeline).}
    \label{fig:pp_colectoresxramalxnodo_col1}
\end{figure}

\begin{figure}[t]
    \centering
    \includegraphics[width=\linewidth]{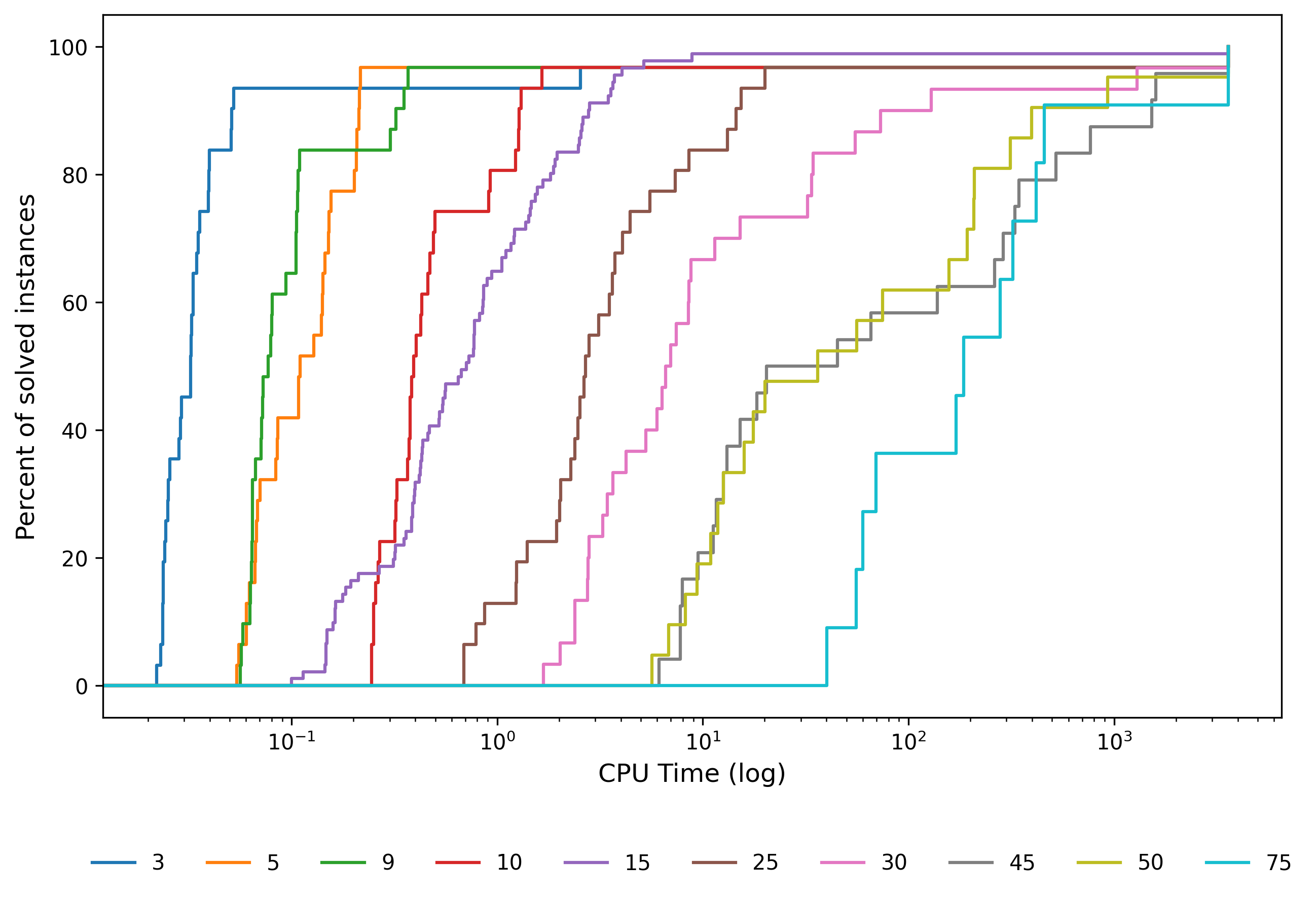}
    \caption{Solution time as a function of branches per pipeline and nodes per branch (two pipelines).}
    \label{fig:pp_colectoresxramalxnodo_col2}
\end{figure}

\begin{figure}[t]
    \centering
    \includegraphics[width=\linewidth]{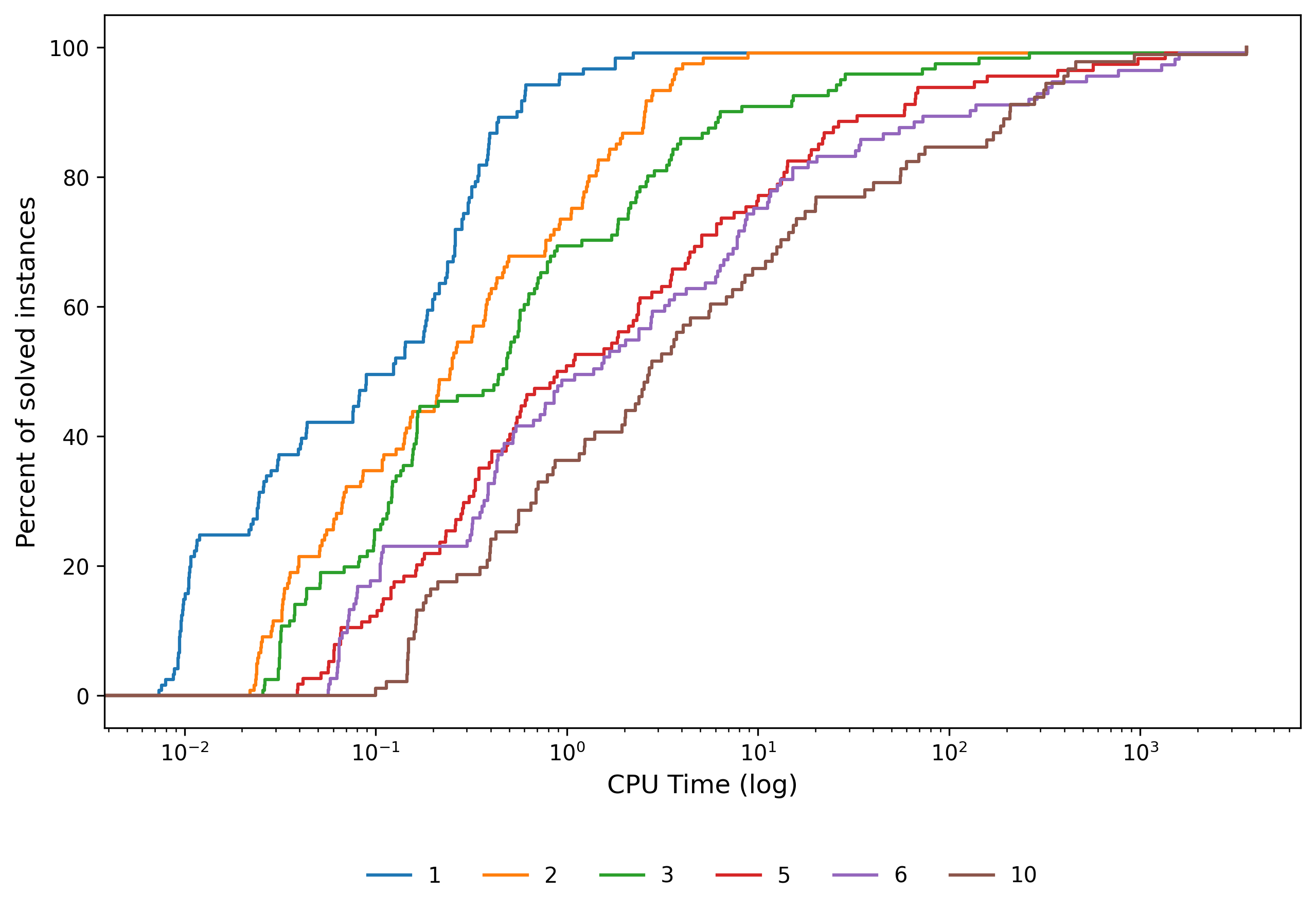}
    \caption{Solution time as a function of the number of pipelines and branches per pipeline (aggregated over nodes).}
    \label{fig:pp_colectoresxramal}
\end{figure}

\begin{figure}[t]
    \centering
    \includegraphics[width=\linewidth]{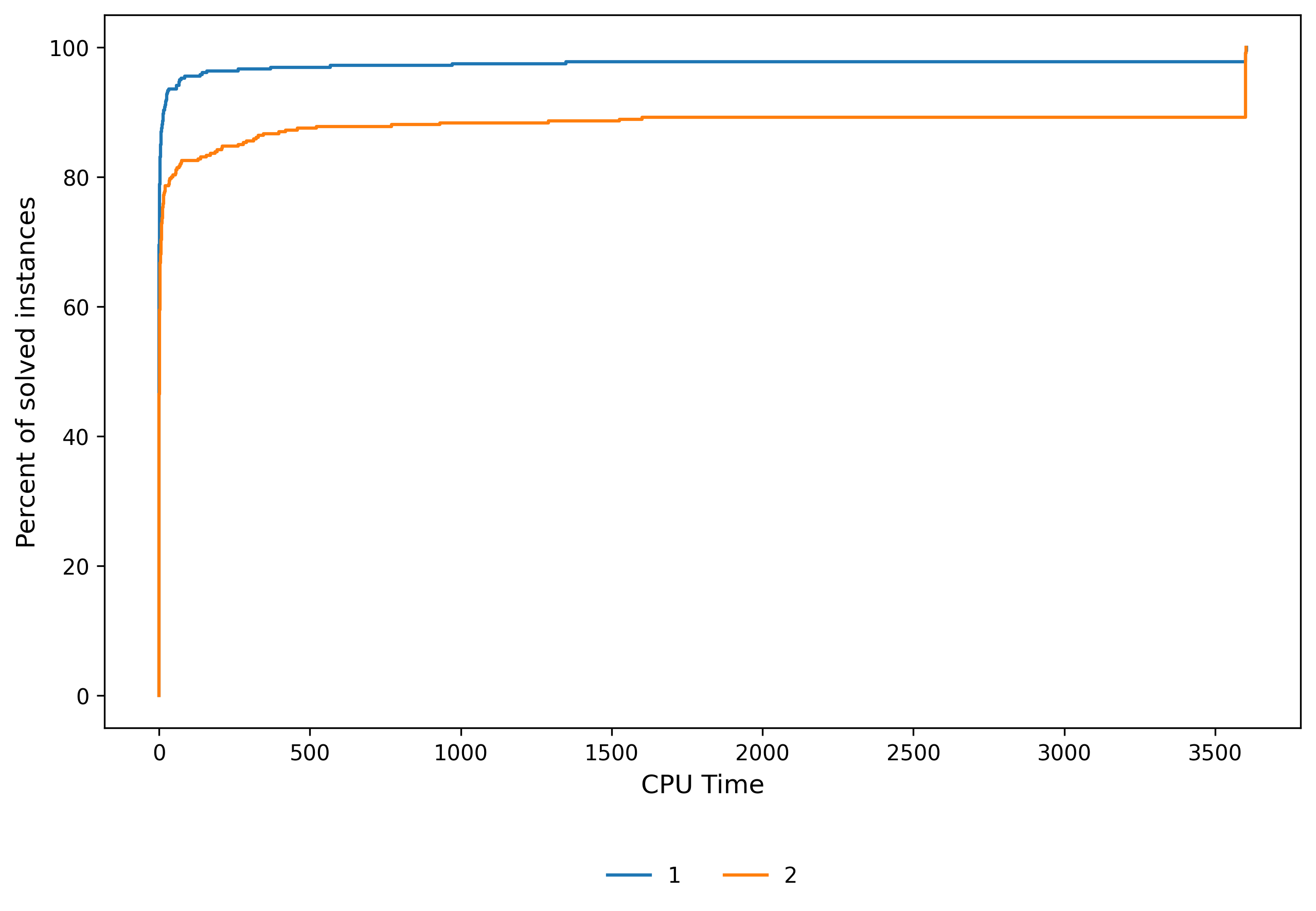}
    \caption{Solution time by number of pipelines (aggregated over branches and nodes).}
    \label{fig:pp_colectores}
\end{figure}

\begin{figure}[t]
    \centering
    \includegraphics[width=\linewidth]{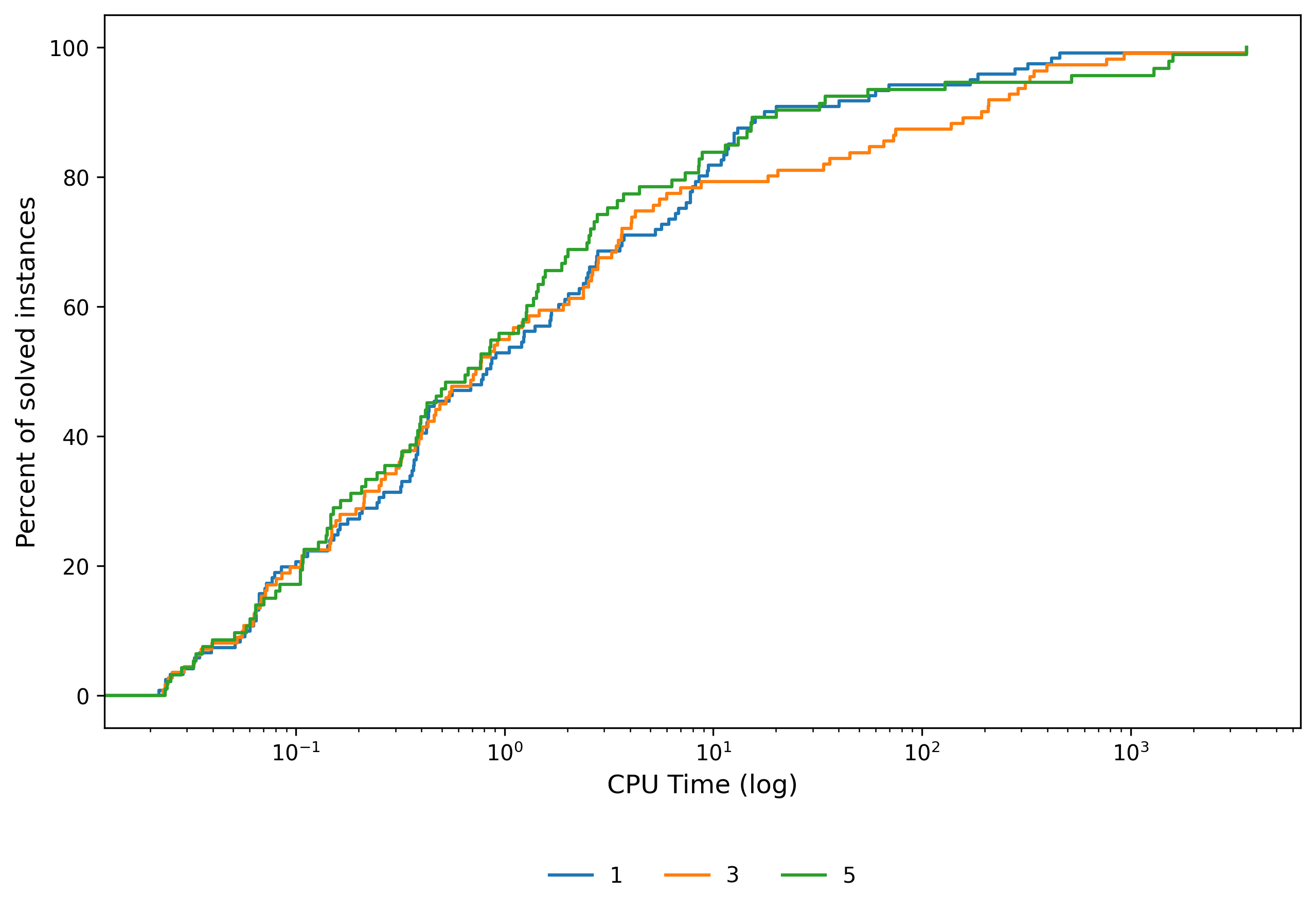}
    \caption{Solution time as a function of safety distance $\Delta$ (aggregated over instance sizes).}
    \label{fig:pp_dist}
\end{figure}

Figures~\ref{fig:pp_colectoresxramalxnodo}--\ref{fig:pp_dist} provide aggregated performance plots that summarize the dependence of solution time on the number of pipelines, branches, nodes, and safety distance.

Overall, the experiments indicate that the proposed discretization--MILP framework can generate feasible wiring layouts efficiently for a broad range of realistic instance sizes, particularly when safety distances are moderate and the number of branches is not excessively large. For highly congested settings (large $\#b$, large $\#n$, and large $\Delta$), computational effort increases sharply and some instances become infeasible within the discretized domain or cannot be solved within the imposed time limit. These findings provide useful guidance for engineering practice, suggesting when exact optimization is practical and when additional decomposition strategies or adaptive discretization may be required.

\section{Industrial Case Study}\label{sec:case}

To demonstrate the practical applicability of the proposed optimization framework, we consider a real engineering scenario provided by our industrial partner, Ghenova, a naval engineering company. The case study concerns the routing of multiple service lines inside a ship cabin, where spatial congestion, structural constraints, and safety regulations make manual design particularly challenging.

The design domain corresponds to a three-dimensional cabin compartment modeled as a rectangular volume of dimensions $5732 \times 2836 \times 2013$ units (in the $X$, $Y$, and $Z$ directions). This volume represents a realistic shipboard installation space in which multiple service pipelines must coexist.

Ten main pipelines (numbered 1--10) are pre-installed in the cabin according to preliminary engineering specifications. These pipelines act as supply backbones from which secondary branches must be routed toward terminal devices located on the boundaries of the compartment.

Five structural walls are present inside the cabin. Each wall consists of a 10-unit-thick metal sheet with eight predefined openings. These walls introduce routing constraints by limiting feasible corridors while allowing controlled passage through designated openings. All routed elements must avoid solid wall regions and may only cross through the permitted openings.

Pipes are modeled as cylindrical elements with a radius of 50 units, which must be respected when enforcing clearance and separation constraints.

\begin{figure}[h]
\begin{center} 
\includegraphics[width=0.4\textwidth, trim=16cm 2cm 14.5cm 2.5cm, clip]{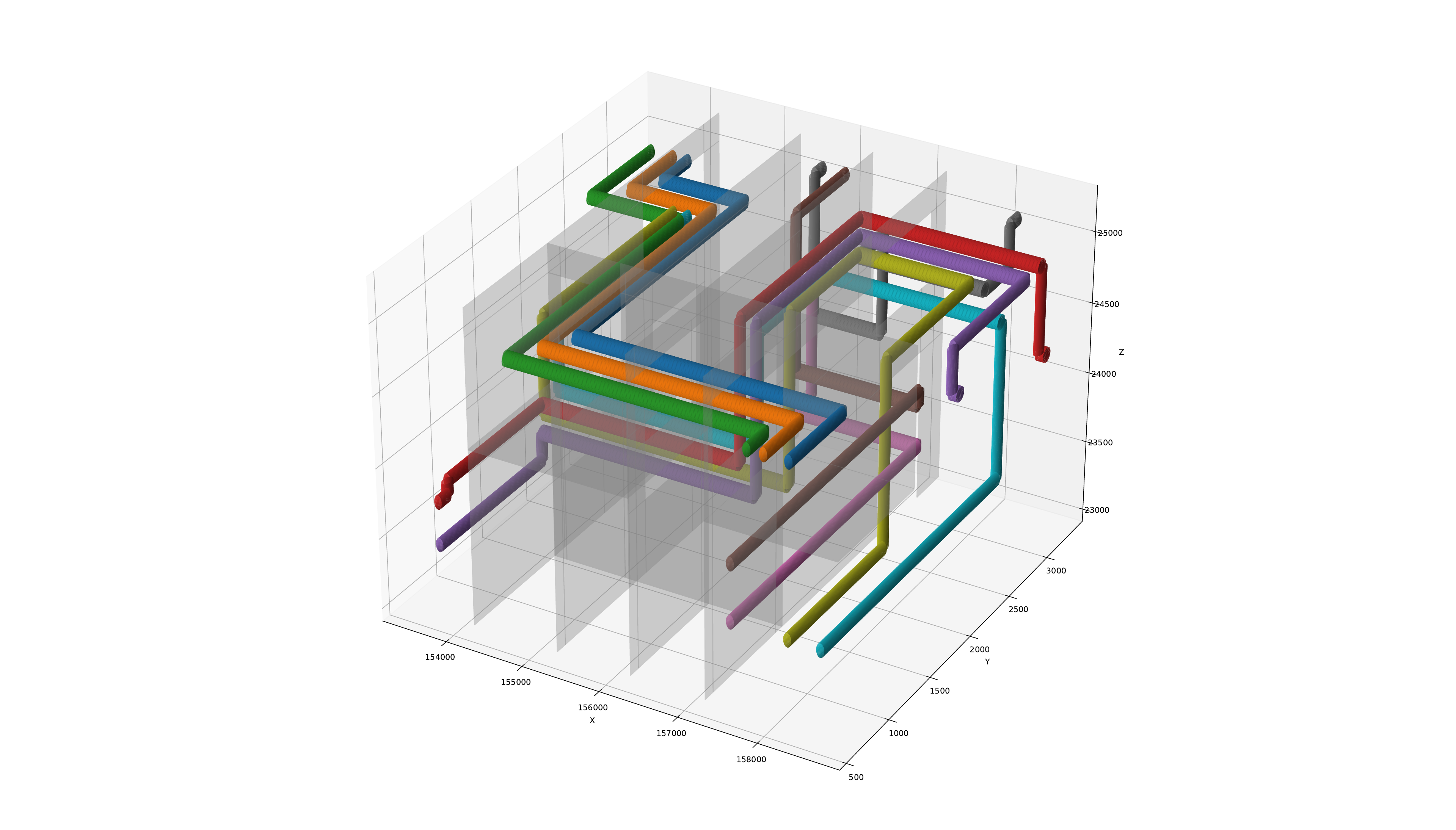}
\includegraphics[width=0.59\textwidth, trim=0.5cm 1cm 0cm 0cm, clip]{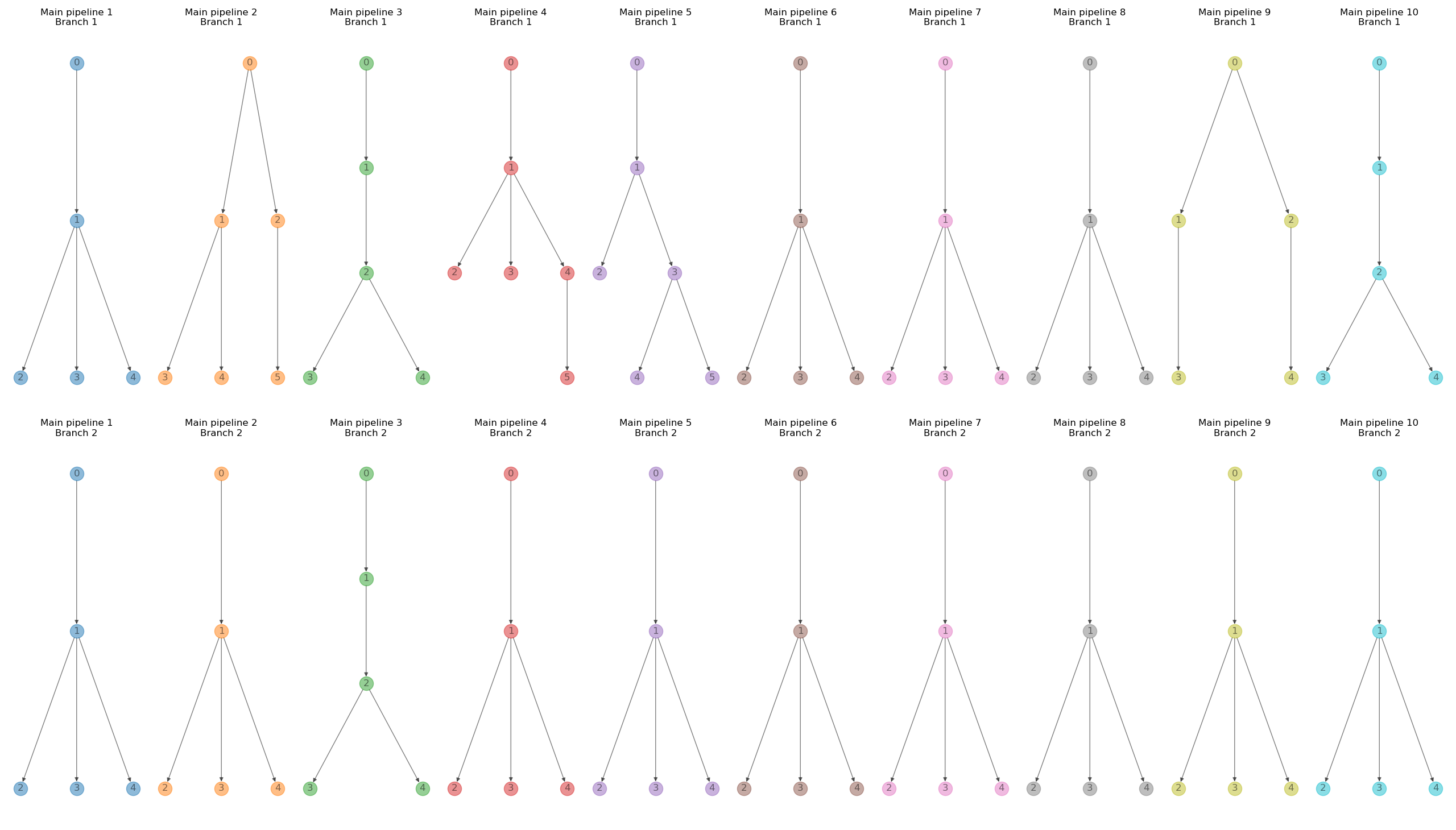}
\end{center}
\caption{Left: three-dimensional cabin layout. Right: predefined branching topologies associated with the main pipelines.} 
\label{fig:entrada}
\end{figure}

Figure~\ref{fig:entrada} illustrates the empty cabin and the required branching topologies.

Each main pipeline is associated with two predefined hierarchical routing schemes that determine how secondary branches must distribute flow to their corresponding terminal devices. These endpoints are located on the boundaries of the cabin.

Every branch must include a valve positioned within a designated cubic region of size $350 \times 350 \times 350$. The exact valve placement within this region is a decision variable of the model. 

To ensure safe installation and maintenance accessibility, a minimum separation distance of 100 units must be maintained between branches. This constraint significantly restricts feasible routing corridors in the already congested environment.

\begin{figure}[h]
\begin{center}
\includegraphics[width=0.49\textwidth, trim=16cm 2cm 14.5cm 2.5cm, clip]{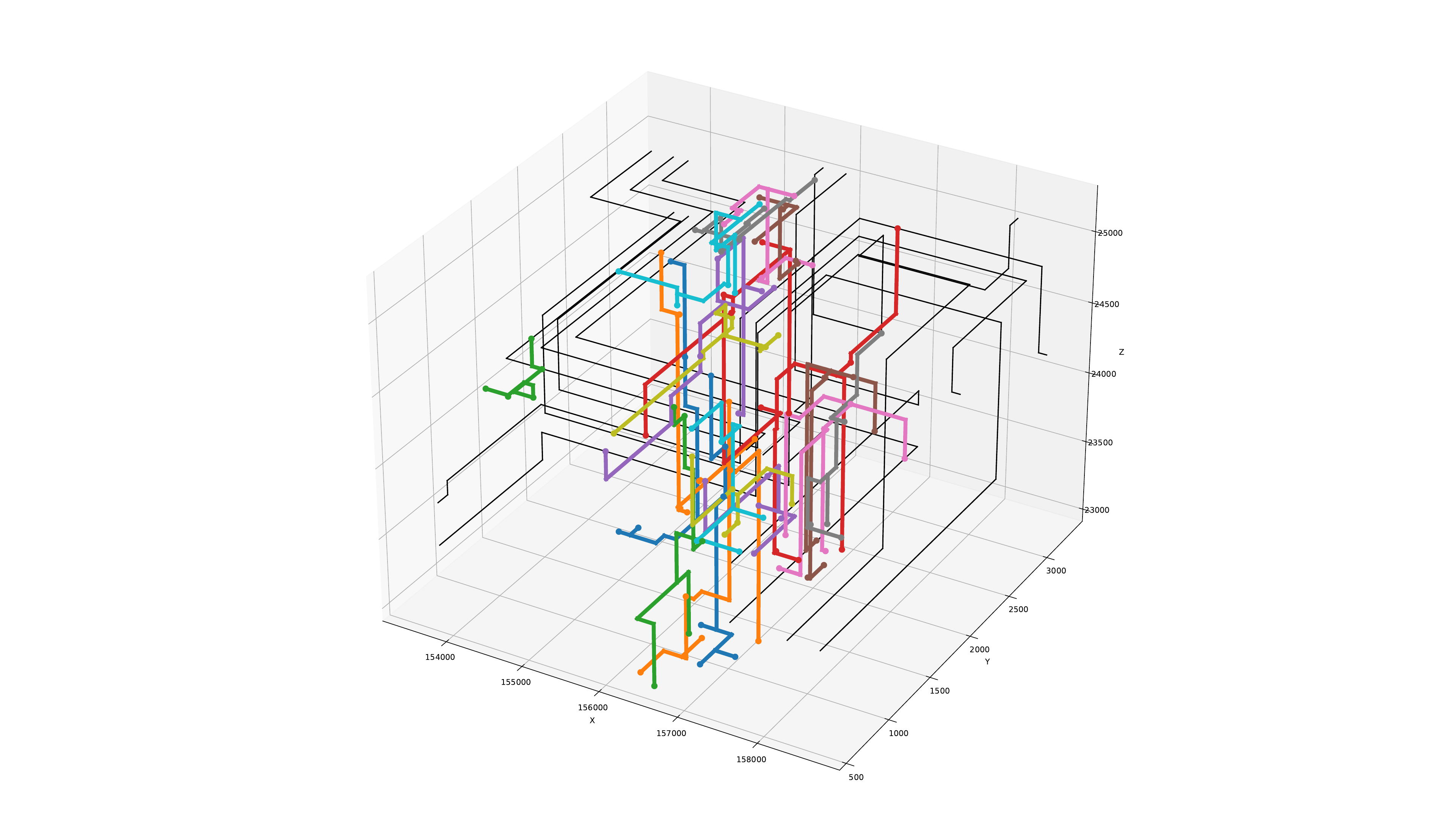}
\includegraphics[width=0.49\textwidth, trim=16cm 2cm 14.5cm 2.5cm, clip]{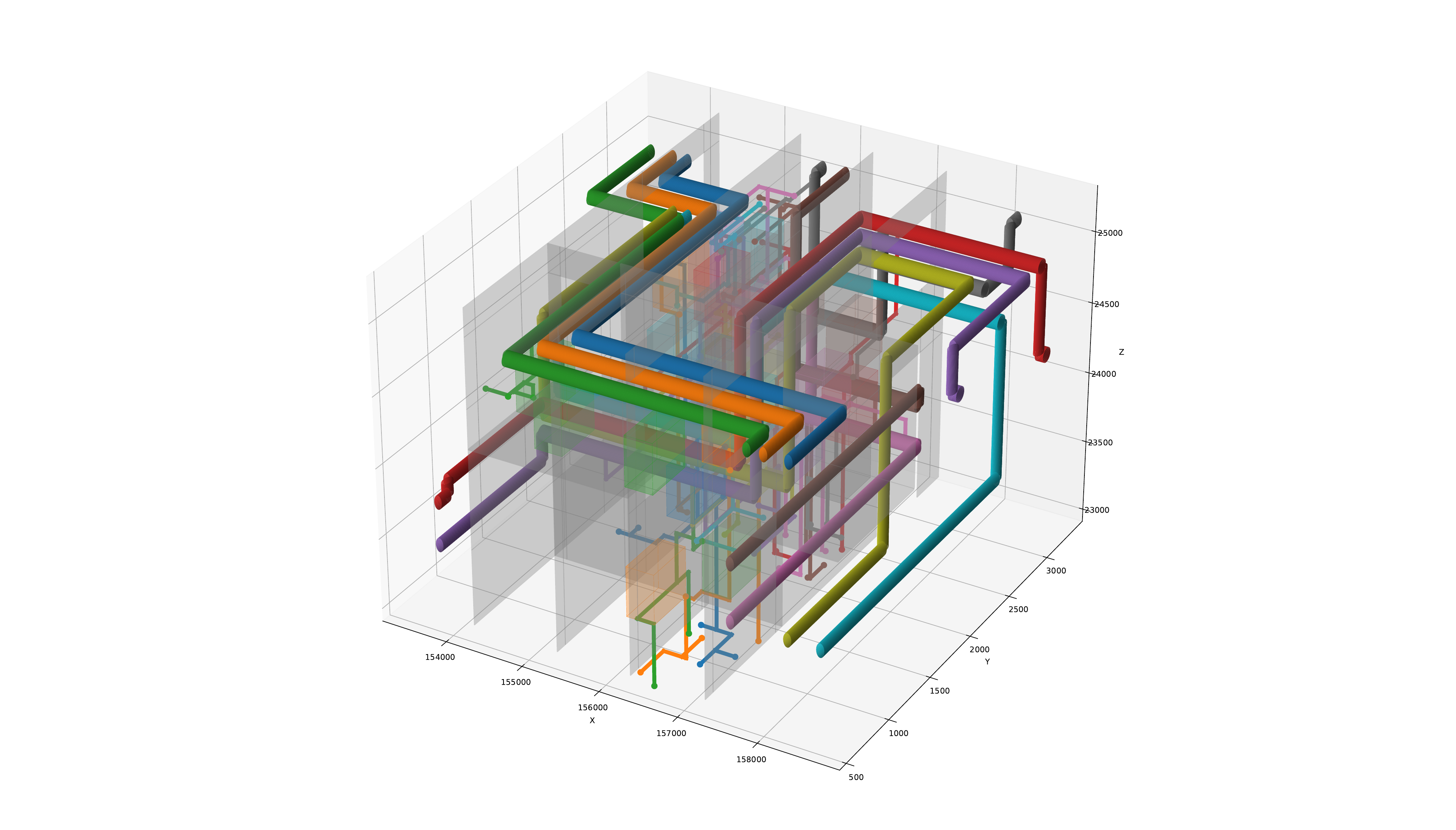}
\end{center}
\caption{Optimized routing layout. Left: main pipelines (black) and routed branches (colored). Right: full three-dimensional configuration including valves and structural walls.} 
\label{fig:salida}
\end{figure}

The left side of Figure~\ref{fig:salida} shows the optimized routing layout obtained by solving the proposed MILP model. Main pipelines are depicted in black, while colored segments represent the routed branches associated with each pipeline. The right-hand panel provides a complete three-dimensional visualization, including valve regions, structural walls, and routed connections.

The instance was solved in 382.88 seconds. The resulting MILP model contained 65,473 binary and flow variables and 100,687 base constraints. During branch-and-bound, 1,106 additional lazy safety constraints were dynamically added, yielding a total of 101,793 constraints. The final objective value (total routed length) was 71,206.0 units.

The case study highlights several practical aspects:

\begin{itemize}
\item The optimization model successfully handles dense spatial configurations with multiple pipelines and structural barriers.
\item Safety distance enforcement substantially shapes routing decisions, often forcing branches to detour through alternative corridors.
\item Valve placement flexibility within admissible regions plays a critical role in enabling feasible layouts.
\item The solution time (under 7 minutes) is compatible with engineering design workflows, enabling the approach to be used as a decision-support tool during layout refinement.
\end{itemize}

Overall, this case study demonstrates that the proposed framework can generate feasible, regulation-compliant wiring layouts for realistic naval engineering environments, providing a systematic alternative to manual routing procedures.

\section{Conclusions}

In this paper, we present an optimization-based framework for the automated design of wiring and pipeline layouts in constrained three-dimensional industrial environments. The proposed approach addresses the Wiring Diagram Problem by integrating hierarchical topology, spatial placement decisions, routing feasibility, and engineering safety requirements within a unified mathematical model.

The proposed methodology combines a structured graph-based discretization of the design space with a mixed-integer linear programming formulation. This integration enables the simultaneous determination of valve placement and branch routing while strictly enforcing technical constraints such as obstacle avoidance, minimum separation distances, and constructibility requirements. By exploiting the hierarchical structure of wiring diagrams, the discretization strategy significantly reduces problem dimensionality and makes exact optimization computationally tractable for realistic engineering instances.

Computational experiments on synthetic configurations demonstrate the scalability of the approach across a wide range of layout densities and safety requirements. The vast majority of instances were solved within a one-hour time limit, with the safety distance parameter identified as the primary driver of computational difficulty for the largest configurations.  The industrial case study, based on a naval cabin environment, further illustrates the practical applicability of the model in handling multiple pipelines, structural walls, valve placement regions, and stringent clearance constraints. The obtained solution times are compatible with engineering design workflows, supporting the use of the framework as a decision-support tool during layout planning and refinement.

Overall, the proposed framework provides a systematic and reproducible alternative to manual routing procedures in complex industrial installations. It enables designers to generate regulation-compliant layouts while minimizing installation length and preserving constructibility.

Future research directions include the development of adaptive discretization strategies to handle highly congested environments more efficiently and the integration of decomposition or heuristic acceleration techniques for very large-scale configurations. Further industrial deployment and validation in additional engineering sectors would also support broader adoption of the methodology.

\section*{Data Availability}
The synthetic instances and source code used in the computational experiments are publicly available at \url{https://github.com/anticiclon/wiring-diagram-problem}. The industrial case study data from Ghenova cannot be shared due to confidentiality agreements.

\section*{Acknowledgements}

The authors acknowledge financial support by  grants PID2020-114594GB-C21, 
PID2022-139219OB-I00, PID2024-156594NB-C21, and 
RED2022-134149-T  (Thematic Network on Location Science and Related Problems) funded by MICIU/AEI/10.13039/501100011033; FEDER+Junta de Andalucía projects C‐EXP‐139‐UGR23, and AT 21\_00032; SOL2024-31596
and SOL2024-31708 funded by US;  the IMAG-Maria de Maeztu grant CEX2020-001105-M / AEI / 10.13039 / 501100011033; and the IMUS--Maria de Maeztu grant CEX2024-001517-M.

\end{document}